\documentclass[a4paper,12pt]{article}
\usepackage{amsmath,amstext,amssymb,mathptmx,fullpage,graphicx,color}
\DeclareMathOperator{\KP}{\mathrm{K}}
\let\le=\leqslant
\let\ge=\geqslant
\newcommand{\qed}{~\raisebox{-0.5ex}{\hbox{\scriptsize$\Box$}}}

\newtheorem{theorem}{Theorem}

\newtheorem{definition}[theorem]{Definition}
\newtheorem{lemma}[theorem]{Lemma}

\begin{document}

\title{Random semicomputable reals revisited}
\author{Laurent~Bienvenu\thanks{%
     LIAFA, CNRS \& Universit\'e Paris Diderot, Paris 7, Case 7014, 75205
     Paris Cedex 13, France,
     e-mail: \texttt{Laurent dot Bienvenu at liafa dot jussieu dot fr}}, 
     Alexander~Shen\thanks{LIRMM, CNRS \& Universit\'e Montpellier 2, 161 rue Ada, 34095 Montpellier Cedex 5, France, on leave from IITP RAS, Bolshoy Karetny, 19, Moscow. e-mail: \texttt{alexander dot shen at lirmm dot fr, sasha dot shen at gmail dot com}. Supported by NAFIT ANR-08-EMER-008-01, RFBR 0901-00709-a grants.}}
\date{October 23, 2011}
\maketitle

\hbox to \hsize{\hss\emph{To Cristian Calude on the occasion on his 60th birthday}}

\begin{abstract}
The aim of this expository paper is to present a nice series of
results, obtained in the papers of Chaitin~\cite{Chaitin1976},
Solovay~\cite{Solovay1975}, Calude et al.~\cite{CaludeHKW1998},
Ku$\mathrm{\check{c}}$era and Slaman~\cite{KuceraS2001}. This joint effort led to a
full characterization of lower semicomputable random reals, both
as those that can be expressed as a ``Chaitin Omega'' and those
that are maximal for the Solovay reducibility. The original proofs
were somewhat involved; in this paper, we present these results
in an elementary way, in particular requiring only basic knowledge
of algorithmic randomness. We add also several simple observations 
relating lower semicomputable random reals and busy beaver functions.
\end{abstract}

\section{Lower semicomputable reals and the
         $\preceq_1$-relation}

Recall that a real number $\alpha$ is \emph{computable} if there is a computable sequence of rationals $a_n$ that converges to $\alpha$ computably: for a given $\varepsilon>0$ one may compute $N$ such that $|a_n-\alpha|\le \varepsilon$ for all~$n>N$. (One can assume without loss of generality that the $a_n$ are increasing.) 

A weaker property is lower semicomputability. A real number $\alpha$ is \emph{lower semicomputable} if it is a limit of a computable increasing sequence of rational numbers. Such a sequence is called \emph{approximation of $\alpha$ from below} in the sequel.

Equivalent definition: $\alpha$ is lower semicomputable if the set of all rational numbers less
than $\alpha$ is enumerable. One more reformulation: if $\alpha=\sum_{i\ge 0} d_i$ where
$d_i$ is computable series of rational numbers, and all $d_i$ with $i>0$ are non-negative.
(We let $d_0$ be negative, since lower semicomputable $\alpha$ can be negative.)

It is easy to see that $\alpha$ is computable if and only if $\alpha$ and $-\alpha$ are lower semicomputable. There exist lower semicomputable but non-computable reals.
Corresponding sequences of rational numbers have non-computable
convergence. (Recall that convergence of a sequence $a_i$ to some $\alpha$
means that for every rational $\varepsilon>0$ there exist some integer 
$N$ such that $|a_i-\alpha|<\varepsilon$ as soon as $i>N$.
Noncomputable convergence means that there is no algorithm that produces
some $N$ with this property given $\varepsilon$.) 

We want to classify computable sequences according to their
convergence speed and formalize the intuitive idea ``one
sequence converges better (i.e., not worse) than the other one''.

\begin{definition}
Let $a_i \to \alpha$ and $b_j \to \beta$ be two computable
strictly increasing sequences converging to lower semicomputable
$\alpha$ and $\beta$ \textup(approximations of $\alpha$ and $\beta$ from below\textup).
We say that 
$a_n\to\alpha$ converges ``better'' \textup(not
worse\textup) than $b_n\to\beta$ if there
exists a total computable function $h$ such that
        $$
\alpha - a_{h(i)} \le \beta - b_i
        $$
for every $i$.
\end{definition}

In other terms, we require that for each term of the second
sequence one may algorithmically find a term of the first one
that approaches the limit as close as the given term of the
second sequence. Note that this relation is transitive (take the
composition of two reducing functions).

In fact, the choice of specific sequences that approximate
$\alpha$ and $\beta$ is irrelevant: \emph{any two increasing
computable sequences of rational numbers that have the same limit, 
are equivalent
with respect to this quasi-ordering}. Indeed, we can just wait
to get a term of a second sequence that exceeds a given term of
the first one. We can thus set the following definition.

\begin{definition}
Let $\alpha$ and $\beta$ be two lower semicomputable reals, and
let $(a_n)$, $(b_n)$ be approximations of $\alpha$ and $\beta$
respectively. If $(a_n)$ converges better than~$(b_n)$, we write
$\alpha \preceq_1 \beta$ \textup(by the above paragraph, this
does not depend on the particular approximations we
chose\textup).
\end{definition}

This definition can be reformulated in different ways. First, we
can eliminate sequences from the defintion and say that
$\alpha\preceq_1\beta$ if there exists a partial computable
function $\varphi$ defined on all rational numbers $r<\beta$
such that
        $$
\varphi(r)<\alpha \text{ \ and \ } \alpha-\varphi(r)\le \beta-r
        $$
for all of them. Below, we refer to $\varphi$ as the
\emph{reduction function}.

The following lemma is yet another characterization of the
order~(perhaps less intuitive but useful).

\begin{lemma}\label{lem:preceq1}
$\alpha\preceq_1\beta$ if and only if $\beta-\alpha$ is lower
semicomputable \textup(or said otherwise, if and only if
$\beta=\alpha+\rho$ for some lower semicomputable
real~$\rho$\textup).
\end{lemma}

\textbf{Proof}. To show the equivalence, note first that
\emph{for every two lower semicomputable reals~$\alpha$ and~$\rho$ 
we have $\alpha\preceq_1\alpha+\rho$}. Indeed, consider
approximations $(a_n)$ to $\alpha$, $(r_n)$ to~$\rho$. Now,
given a rational $s < \alpha+\rho$, we wait for a stage~$n$ such
that $a_n+r_n>s$. Setting $\varphi(s)=a_n$, it is easy to check
that $\varphi$ is a suitable reduction function witnessing
$\alpha \preceq_1 \alpha+\rho$.

It remains to prove the reverse implication: \emph{if
$\alpha\preceq_1\beta$ then $\rho=\beta-\alpha$ is lower
semicomputable}. Indeed, if $(b_n)$ is a computable
approximation (from below) of $\beta$ and $\varphi$ is the
reduction function that witnesses $\alpha\preceq_1\beta$, then
all terms $b_n - \varphi(b_n)$ are less than or equal to
$\beta-\alpha$ and converge to $\beta-\alpha$. (The sequence 
$b_n-\varphi(b_n)$ may
not be increasing, but still its limit is lower semicomputable,
since all its terms do not exceed the limit, and we may replace
$n$th term by the maximum of the first $n$ terms.)\qed

\smallskip
A special case of this lemma: let $\sum u_i$ and $\sum v_i$ be computable
series with non-negative rational 
terms (for $i>0$; terms $u_0$ and $v_0$ are starting points
and may be negative) that converge to (lower semicomputable) $\alpha$ and $\beta$. 
If $u_i\le v_i$ for all $i>0$, then $\alpha\preceq_1\beta$, since 
$\beta-\alpha=\sum_i (v_i-u_i)$ is lower semicomputable.

The reverse statement is also true: if $\alpha\preceq_1\beta$, one can find computable series $\sum u_i=\alpha$ and $\sum v_i=\beta$ with these properties ($0\le u_i\le v_i$ for $i>0$). Indeed, $\beta=\alpha+\rho$ for lower semicomputable $\rho$; take $\alpha=\sum u_i$ and $\rho=\sum r_i$ and let $v_i=u_i+r_i$.

In fact, a stronger statement is also true; each of 
the series can be chosen in an arbitrary way. We have already
seen how to choose $v_i$ when $u_i$ are given. The other direction:
assume that $\alpha\preceq_1\beta=\sum v_i$ for some $v_i\ge 0$.  
We need a decomposition
$\alpha=\sum u_i$ where $u_i\ge 0$ and $u_i\le v_i$ for $i>0$. 
Indeed, we can construct $u_i$ sequentially
using the following invariant: the current approximation 
$A=\sum_{j<i} u_j$ to $\alpha$ should be below $\alpha$ and at 
least as close (to $\alpha$)
as the current approximation $B=\sum_{j<i}v_j$ (to $\beta$). Initially we choose $u_0$ applying reduction function to $v_0$. When the current approximation becomes $B'=B+v_i$, we apply reduction function to get $A'$ which is at least as close to $\alpha$ as $B'$ is to $\beta$. Then there are several cases: 

(1)~if $A'<A$, we let $u_i=0$, and the next approximation is $A$ (it is close enough by assumption);

(2)~if $A\le A'\le A+v_i$, we let $u_i=A'-A$; the condition guarantees that $u_i\le v_i$;

(3)~finally, if $A'>A+v_i$, we let $u_i=v_i$ (the invariant remains valid since the distances to $\alpha$ and $\beta$ are decreased by the same amount). 

\section{The Solovay reducibility and complete reals}

Let $\alpha$ be a lower semicomputable but not computable real.
By the results of the previous section, one has
      \[
\alpha\preceq_1 2\alpha\preceq_1 3\alpha \preceq_1 \ldots
      \]
because for all~$k$ the difference $(k+1)\alpha-k\alpha=\alpha$\label{pg:2alpha}
is lower semicomputable (so Lemma~\ref{lem:preceq1} applies).
The reverse relations are not true, because
$k\alpha-(k+1)\alpha=-\alpha$ is not lower semicomputable (if it
were, then $\alpha$ would be computable).

One may argue that this relation is therefore a bit too sharp.
For example, $\alpha$ and $2\alpha$ have essentially the same
binary expansion (just shifted by one position), so one may want
$\alpha$ and $2\alpha$ to be equivalent. In other words, one may
look for a less fine-grained relation. A natural candidate for
this is \emph{Solovay reducibility}.

\begin{definition}[Solovay reducibility]
We say that $\alpha\preceq\beta$ if $\alpha\preceq_1 c\beta$ for
some positive integer $c>0$.
\end{definition}

(A convenient notation: we say, for some positive rational~$c$, that $\alpha\preceq_c \beta$ if
$\alpha\preceq_1 c\beta$. Then $\alpha\preceq\beta$ if
$\alpha\preceq_c\beta$ for some~$c$.)

Like for lower semicomputable semimeasures in algorithmic information theory (see, e.g.,~\cite{shen2000a}), one can easily prove
the existence of maximal elements~\cite{Solovay1975}.

\begin{theorem}
There exists a $\preceq$-biggest lower semicomputable real.
\end{theorem}

\textbf{Proof}.
Indeed, we can enumerate all lower semicomputable reals
$\alpha_i$ in $[0,1]$ and then take their sum $\alpha=\sum
w_i\alpha_i$ with computable positive weights $w_i$ such that
$\sum w_i$ converges. This $\alpha$ can be represented as
$w_i\alpha_i$ plus some lower semicomputable real, so $\alpha_i
\preceq_1 (1/w_i)\alpha$.\qed

\medskip

The biggest elements for the $\preceq$-preorder are also called
\emph{(Solovay) complete} lower semicomputable reals. They have
an alternative description~\cite{Solovay1975,CaludeHKW1998}:

\begin{theorem}
        \label{omega-complete}
Complete semicomputable reals in $[0,1]$ are sums of universal semimeasures
on $\mathbb{N}$ and vice versa.
\end{theorem}

Recall (see~\cite{shen2000a} for details) that lower semicomputable semimeasures on~$\mathbb{N}$ are lower semicomputable functions $m\colon\mathbb{N}\to\mathbb{R}$ with non-negative values such that $\sum_i m(i)\le 1$. (For a function $m$ lower semicomputability means that $m(i)$ is lower semicomputable uniformly in $i$: there is an algorithm that gets $i$ as input and produces an increasing sequence of rationals that converges to $m(i)$.) Universal semimeasures are the maximal (up to a constant factor) lower semicomputable semimeasures.

\textbf{Proof}. Any lower semicomputable real $\alpha$ is a sum
of a computable series of rationals; this series (up to a
constant factor that does not matter due to the definition of the
Solovay reducibility) is bounded by a universal semimeasure. The
difference between the upper bound and the series itself is a
lower semicomputable semimeasure, and therefore $\alpha$ is
reducible to the sum of the universal semimeasure.

We have shown that sums of universal semimeasures are complete.
On the other hand, let $\alpha$ be a Solovay complete real in $[0,1]$.
We need to show that $\alpha$ is a sum of some universal semimeasure.  
Let us start with arbitrary universal semimeasure $m(i)$. The sum
$\sum m(i)$ is lower semicomputable and therefore $\sum m(i)\preceq_1
c\alpha$, so $\alpha=\sum m(i)/c + \tau$ for some integer $c>0$ and
some lower semicomputable $\tau$. Dividing $m$ by $c$ and then adding
$\tau$ to one of the values, we get a universal semimeasure
with sum $\alpha$.\qed

Chaitin denoted the sum of a universal semimeasure by $\Omega$.
Since there is no such thing as \emph{the} universal
semimeasure, it is better to speak about \emph{$\Omega$-reals} defined
as sums of universal semimeasures. We have shown therefore that
the class of $\Omega$-reals coincides with the class of Solovay complete
lower semicomputable reals in $[0,1]$.

It turns out that this class has one more characterization~\cite{Chaitin1976,CaludeHKW1998,KuceraS2001}:  
\begin{theorem}
	A lower semicomputable real is complete if and only if it is
	Martin-L\"of random.
\end{theorem}

(See, e.g.,~\cite{shen2000a} for the definition of Martin-L\"of randomness.)
We provide the proof of this result below, starting with one direction in the next section~\ref{sec:complete-random} and finishing the other direction in section~\ref{sec:random-complete}. 

\section{Complete lower semicomputable reals are random}\label{sec:complete-random}

The fact that lower semicomputable reals are random, is Chaitin's theorem (randomness of $\Omega$). It is usually proved by using complexity characterization of
randomness. However, there is a direct argument that does not
involves complexity (it is in the footnote in Levin's
``Forbidden information'' paper~\cite{LevinTA1}; this footnote
compressed the most important facts about lower semicomputable
random reals into few lines!).

First, we prove that \emph{there exists a lower semicomputable
random real}. For that we consider an effectively open set $U$
of measure less than (say) $1/2$ that covers all non-random
reals in $[0,1]$. (The definition of Martin-L\"of randomness
guarantees that for every $\varepsilon>0$ one can find an
effectively open set that has measure less than~$\varepsilon$
and covers all non-random reals. We need only one such set for
some $\varepsilon<1$, say, $\varepsilon=1/2$.) Then take the
minimal element $\alpha$ in a closed set $[0,1]\setminus U$.
This number is random (by definition) and lower semicomputable:
compactness implies that any segment $[0,r]$ with
rational $r<\alpha$ is covered by finitely many intervals of~$U$
and thus all such $r$'s can be enumerated.

Second, we prove that \emph{randomness is upward-closed}: if
$\alpha\preceq\beta$ and $\alpha$ is random, then $\beta$ is
random. We may assume without loss of generality that
$\alpha\preceq_1\beta$ (randomness does not change if we
multiply a real by a rational factor).

So let $b_i\to\beta$ be a computable increasing sequence of
rational numbers that converges to $\beta$. Assume that somebody
gives us (in parallel with $b_i$) a sequence of rational intervals and guarantees that
one of them covers $\beta$. How to transform it into a sequence
of intervals that covers $\alpha$ (i.e., one of the intervals
covers~$\alpha$) and has the same (or smaller) total length?
If an interval appears that is entirely on the
left of the current approximation $b_i$, it can be ignored
(since it cannot cover $\beta$ anyway). If the interval is entirely on
the right of $b_i$, it can be postponed until the current
approximation $b_j$ enters it (this may happen or not, in the
latter case the interval does not cover $\beta$). If the
interval contains $b_i$, we can convert it into the interval of
the same length that starts at $a_j$, where $a_j$ is a rational
approximation to $\alpha$ that has the same or better precision
as $b_i$ (as an approximation to~$\beta$): if $\beta$ is in
the original interval, $\alpha$ is in the converted interval.

So randomness is upward-closed and therefore complete lower
semicomputable reals are random.

\textbf{Remark}. The second part can be reformulated: if $\alpha$ and $\beta$ are lower semicomputable reals and at least one of them is random, then the sum $\alpha+\beta$ is random, too. The reverse is also true: if both $\alpha$ and $\beta$ are non-random, then $\alpha+\beta$ is not random. (We will see later different proofs of this statement.)

\section{Randomness and prediction game}

Before proving the reverse implication, let us make a digression and look more closely at the last argument. Consider the following game: an observer watches an increasing sequence of rationals (given one by one) and from time to time makes predictions of the following type: ``the sequence will never increase by more than $\delta$'' (compared to its current value). Here $\delta$ is some non-negative rational. The observer wins this game if (1)~one of the predictions remains true forever; (2)~the sum of all numbers $\delta$ used in the predictions is small (less that some rational $\varepsilon>0$ which is given to the observer in advance).

It is not required that at any moment a valid prediction exists, though one could guarantee this by making predictions with zero or very small (and decreasing fast) $\delta$ at each step. Note also that every prediction can be safely postponed, so we may assume that the next prediction is made only if the previous one becomes invalid.  Then at any moment there is only one valid prediction.

\begin{theorem}
Let $a_i$ be a computable increasing sequence of rational numbers that converges to some \textup(lower semicomputable\textup) real~$\alpha$. The observer has a computable winning strategy in the game if and only if $\alpha$ is not random.
\end{theorem}

\textbf{Proof}. A computable winning strategy gives us a computable sequence of prediction intervals of small total measure and guarantees that one of these (closed) intervals contains~$\alpha$. On the other hand, having a sequence of intervals that covers $\alpha$ and has small total measure, we may use it for predictions.  To make the prediction, we wait until the current approximation $a_i$ gets into the already discovered part of the cover (this will happen since the limit is covered). Then for our prediction we use the maximal $\delta$ such that $(a_i,a_i+\delta)$ is covered completely at the moment, and then wait until this prediction becomes invalid. Then the same procedure is used again. At some point $\alpha$ is covered by some interval in the sequence and the current approximation enters this interval; the prediction made after this moment will remain valid forever. The total length of all prediction interval is bounded by the measure of the cover (the prediction intervals are disjoint and all are covered).\qed

A reformulation of the same observation that does not use game terminology:

\begin{theorem}
	\label{game-criterion}
Let $a_i$ be a computable increasing sequence of rational numbers that converges to~$\alpha$. The number $\alpha$ is non-random if and only if for every rational $\varepsilon>0$ one can effectively find a computable sequence $h_0,h_1,\ldots$
of non-negative rational numbers such that $\sum_i h_i<\varepsilon$ and $\alpha\le a_i+h_i$ for some $i$.
\end{theorem}

(Here the predictions $h_i$ are made on every step; it does not matter
since we may use zeros.)

There is a Solovay criterion of randomness (a constructive version of
Borel--Cantelli lemma): a real number $\alpha$ is non-random if and only if
there exists a computable sequence of intervals that have finite total
measure and cover $\alpha$ infinitely many times. It can also be
reformulated in the style of our previous theorem:

\begin{theorem}
Let $a_i$ be a computable increasing sequence of rational numbers that converges to~$\alpha$. The number $\alpha$ is non-random if and only if  there exists a computable sequence $h_0,h_1,\ldots$ of non-negative rational numbers such that $\sum_i h_i<\infty$ and $\alpha\le a_i+h_i$ for infinitely many $i$.
\end{theorem}

\textbf{Proof}. If $\alpha$ is non-random, we apply the preceding result for $\varepsilon=1, 1/2, 1/4, 1/8,\ldots$ and then add the resulting sequences (with shifts $0,1,2,\ldots$ to the right). Each of them provides one value of $i$ such that $\alpha\le a_i+h_i$, and these values cannot be bounded due to shifts. On the other hand, if $\alpha\le a_i+h_i$ for infinitely many~$i$, we get a sequence of intervals with finite sum of measures that covers $\alpha$ infinitely many times (technically, we should replace closed intervals by slightly bigger open intervals). It remains to use Solovay's criterion (or recall its proof: the effectively open set of points that are covered with multiplicity $m$ has measure at most $O(1/m)$).\qed

\medskip

The randomness criterion given in this section implies the following observation (which may look strange at first). Consider a sum of a computable series of positive rational numbers. \emph{The randomness of the sum cannot change if all summands are changed by some~$\Theta(1)$-factor}. Indeed, all $h_i$ can be multiplied by a constant.

Now let us prove that \emph{if $\alpha$ and $\beta$ are non-random lower semicomputable reals, their sum 
$\alpha+\beta$ is non-random, too}. (See the discussion in the previous section). The natural idea to prove this is the following:
make predictions in the games for $\alpha$ and $\beta$, and then take their sum as prediction for 
$\alpha+\beta$. But this simple argument does not work. The problem is that the same prediction for $\alpha$ can be combined
with many predictions for $\beta$ and therefore will be counted many times in the sum.

The solution is to make predictions for $\alpha$ and $\beta$ of the same size. Let $a_i$ and $b_i$ 
be computable increasing sequences that converge to $\alpha$ and $\beta$.
Since $\alpha$ and $\beta$ are non-random, they are covered by sequences of intervals that have
small total measure. 
To make a prediction for the sequence $a_i+b_i$ (after the previous prediction became invalid) we wait until the current approximations $a_i$ and $b_i$ become covered by the intervals of those sequences. We take then the maximal $h$ and $k$ such that $(a_i,a_i+h)$ and $(b_i,b_i+k)$ are entirely covered (by the unions of currently appeared intervals). The prediction interval is declared to be
$(a_i+b_i, a_i+b_i+\delta)$ where $\delta = 2\min(h,k)$.

Let us show that one of the predictions will remain valid forever. Indeed, the limit values $\alpha$ and $\beta$ are covered by some intervals. These intervals appear in the sequences at some point and cover $\alpha$ and $\beta$ with some neighborhoods, say, $\sigma$-neighborhoods. If the prediction is made after $a_i$ and $b_i$ enter these neighborhoods, $\delta$ is greater than $2\sigma$ and the prediction is final: $a_i+b_i$ never increases more than by $\delta$.

It remains to estimate the sum of all $\delta$s used during the prediction. It can be done using the following observation: when a prediction interval $(a_i+b_i, a_i+b_i+\delta)$ becomes invalid, this means that either $a_i$ or $b_i$ has increased by $\delta/2$ or more, so the total measure of the cover on the right of $a_i$ and $b_i$ has decreased at least by $\delta/2$. (Here we use that $(a_i,a_i+\delta/2)$ and $(b_i, b_i+\delta/2)$ are covered completely because $\delta/2$ does not exceed both $h$ and $k$: it is important here that we take the minimum.)

\medskip
Let us return to the criterion for randomness provided by Theorem~\ref{game-criterion}. The condition for non-randomness given there can be weakened in two aspects: first, we can replace computable sequence by a semicomputable sequence; second, we can replace $h_i$ by the entire tail $h_i+h_{i+1}+\ldots$ of the corresponding series:

\begin{theorem}
	\label{semicomputable-predictions}
Let $a_i$ be an increasing computable sequence of rational numbers that converges to~$\alpha$. Assume that for every rational $\varepsilon>0$ one can effectively find a lower semicomputable sequence $h_i$ of non-negative reals such that $\sum_i h_i < \varepsilon$ and $\alpha \le a_i +h_i+h_{i+1}+\ldots$ for some $i$. Then $\alpha$ is not random.
\end{theorem}

\textbf{Proof}. Assume that for every $i$ there is a painter who get $h_i$ units of paint and the instruction to paint the line starting at $a_i$, going to the right and skipping the parts already painted by other painters (but making no other gaps). (Since $h_i$ is only semicomputable, the paint is provided incrementally.) The painted zone is an effective union of intervals of total measure $\sum_i h_i$. If $\alpha<a_i+h_i+h_{i+1}+\ldots$, then $\alpha$ is painted since we cannot use $h_i+h_{i+1}+\ldots$ paint starting between $a_i$ and $\alpha$ (recall that all $a_k$ are less than $\alpha$) and not crossing $\alpha$. (In the condition we have $\le$ instead of $<$, but this does not matter since we can increase all $h_i$ to, say, twice their original value.)\qed

\medskip

This result implies one more criterion of randomness for lower semicomputable reals:

\begin{theorem}
Let $\alpha=\sum d_i$ be a computable series of non-negative rational numbers. The number $\alpha$ is non-random if and only if for every $\varepsilon>0$ one can effectively produce an enumerable set $W\subset\mathbb{N}$ of indices such that \textup{(1)}~$\sum_{i\in W} d_i < \varepsilon$
and \textup{(2)}~$W$ is co-finite, i.e., contains all sufficiently large integers.
\end{theorem}

\textbf{Proof}. If $\alpha$ is not random, it can be covered by intervals with arbitrarily small total measure. It remains to consider the set $W$ of all $i$ such that $(d_0+\ldots+d_{i-1}, d_0+\ldots+d_{i-1}+d_i)$ is entirely covered by one of those intervals. In the other direction the statement is a direct consequence of Theorem~\ref{semicomputable-predictions}, just let $a_i=d_0+\ldots+d_{i-1}$ and $h_i=d_i$ for $i\in W$ (and $h_i=0$ for $i\notin W$).\qed

\medskip

This result shows again that the sum of two non-random lower semicomputable reals is not random (take the intersection of two sets $W_1$ and $W_2$ provided by this criterion for each of the reals).

\section{Random lower semicomputable reals are complete}\label{sec:random-complete}

To prove the completeness of random lower semicomputable reals, let us start with the following remark. Consider two lower semicomputable reals
$\alpha$ and $\beta$ presented as  limits of
increasing computable sequences $a_i\to\alpha$ and
$b_i\to\beta$.
Let $h_i=a_{i+1}-a_i$ be the increases in the first sequence.
We may use $h_i$ to construct a strategy for the prediction game against the second sequence in the following way. We shift the interval $[a_1,a_2]$ to get the (closed) interval of the same length that starts at $b_1$. Then we wait until $b_i$ at the right of this interval
appears; let it be $b_{i_1}$. Then shift the interval
$[a_2,a_3]$ to get the interval of the same length that starts
at $b_{i_1}$; let $b_{i_2}$ be the first $b_i$ on the right of
it, etc.

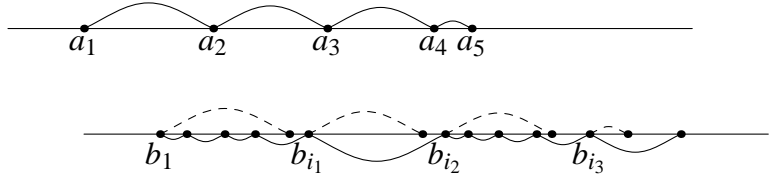
\begin{center}
        \includegraphics[scale=1.0]{lowerrandom-1.mps}
\end{center}

There are two possibilities: either 

(1)~the observer wins in the
prediction game, i.e., some of the shifted
intervals covers the rest of $b_i$ and the next $b_{i_k}$ is
undefined, or
 
(2)~this process continues indefinitely.

In the second case
$\alpha\preceq_1\beta$ since the difference $\beta-\alpha$ is
represented as a sum of a computable series (``holes'' between
neighbor intervals; note that the endpoints of the shifted intervals
also converge to~$\beta$).

One of these two alternatives happens for arbitrary lower semicomputable
reals $\alpha$ and~$\beta$. Now assume that $\beta$ is not Solovay
complete; we need to prove that $\beta$ is not random.
Since $\beta$ is not complete, there exists some $\alpha$ such that
$\alpha\not\preceq\beta$. In particular, $\alpha\not\preceq_1\beta$.
Therefore, for these $\alpha$ and $\beta$ the second
alternative is impossible, and the observer wins. In other terms,
we get a computable sequence of
(closed) intervals that covers $\beta$. Repeating the same
argument for $\alpha/2$, $\alpha/4$,\ldots\ (we know that 
$\alpha/c\not\preceq_1\beta$ for every $c$, since $\alpha\not\preceq\beta$)
we effectively get a
cover of $\beta$ with arbitrary small measure (since the sum of 
all $h_i$ is bounded by a integer
constant even being non-computable),
therefore $\beta$ is not random.

\textbf{Remark}. This argument probably gives some quantitative
connection between randomness deficiency of a random lower semicomputable
real and
another parameter that can be called \emph{completeness
deficiency}. It can be defined as follows: fix some complete
$\alpha$ and for every $\beta$ consider the infimum of all $c$ such
that $\alpha\preceq_1 c\beta$. 

\section{Slow convergence: Solovay functions}

We have seen several results of the following type: the limit of an increasing computable sequence of rationals is random if and only if the convergence is slow. In this section we provide one more result of this type.

Consider a computable converging series $\sum r_i$ of positive rational numbers. Note that $r_i$ is bounded by $O(m_i)$ where $m\colon i\mapsto m_i$ is a universal semimeasure ($m_i$ is also called \emph{a priori probability} of integer $i$). Therefore prefix complexity $\KP(i)=-\log_2 m_i$ is bounded by $-\log_2 r_i + O(1)$ (see, e.g,~\cite{shen2000a}). We say that the series $\sum r_i$ \emph{converges slowly in the Solovay sense}  (has the \emph{Solovay property}) if this bound is tight infinitely often, i.e., if $r_i \ge \varepsilon m_i$ for some $\varepsilon>0$ and for infinitely many $i$. In other terms, the series does \emph{not} converge slowly if $r_i/m_i\to 0$.

In~\cite{BienvenuD2009,HolzlKM2009} the name \emph{Solovay function} was used for a computable bound $S(i)$ for prefix complexity $\KP(i)$ that is tight infinitely often, i.e., $\KP(i)\le S(i)+O(1)$ for every $i$ and $\KP(i)\ge S(i)-c$ for some $c$ and for infinitely many values of $i$. Thus, a computable series $\sum a_i$ of positive rational numbers has the Solovay property  if and only if $i\mapsto -\log_2 a_i$ is a Solovay function~\cite{BienvenuD2009}.

\begin{theorem}
 	\label{Solovay-property-and-randomness}
Let $\alpha=\sum_i r_i$ be a computable converging series of positive rational numbers. The number $\alpha$ is random if and only if this series converges slowly in the Solovay sense.
\end{theorem}

In other terms, the sum is non-random if and only if the ratio $r_i/m_i$ tends to~$0$.

\textbf{Proof}. Assume that $r_i/m_i\to 0$. Then for every $\varepsilon$ we can let $h_i=\varepsilon m_i$ and get a lower semicomputable sequence that satisfies the conditions of Theorem~\ref{semicomputable-predictions}. Therefore $\alpha$ is not random.

We can also prove that $\alpha$ is not complete (thus providing an alternative proof of its non-randomness). Recall the argument used in the proof of Theorem~\ref{omega-complete}: if $r_i \le m_i$, then $\sum r_i\preceq_1\sum m_i$. And if $r_i \le c m_i$, then $\sum r_i \preceq_c \sum m_i$. This remains true if the inequality $r_i \le c m_i$ is true for all sufficiently large $i$. So for a fast (non-Solovay)  converging series and its sum $\alpha$ we have $\alpha \preceq_c \sum m_i$ \textsl{for arbitrarily small~$c$}. If $\alpha$ were complete, we would have also $\sum m_i \preceq_d \alpha$ for some $d$ and therefore $\alpha\preceq_{cd}\alpha$ for some $d$ and all $c>0$. For small enough $c$ we have $cd <1/2$
and therefore $\alpha\preceq_{1/2}\alpha$ i.e., $2 \alpha \preceq_1 \alpha$. Then, as we saw on page~\pageref{pg:2alpha}, $\alpha$ should be computable.

It remains to show the reverse implication. Assuming that $\alpha=\sum r_i$ is not random, we need to prove that $r_i/m_i\to 0$. Consider the interval $[0,\alpha]$ split into intervals of length $r_0,r_1,\ldots$. Given an open cover of $\alpha$ with small measure, we consider those intervals (of length $r_0,r_1,\ldots$, see above) that are
completely covered (endpoints included). They form an enumerable set and the sum of their lengths does not exceed the measure of the cover. 
If the cover has measure $2^{-2n}$ for some $n$, we may multiply the corresponding $r_i$ by $2^n$ and their sum remains at most $2^{-n}$. Note also that for large enough $i$ the $i$th interval  is
covered (since it is close to $\alpha$ and $\alpha$ is covered). So for each $n$ we get a semimeasure $M^n=M^n_0,M^n_1,\ldots$
such that $M^n_i/r_i\ge 2^n$ for sufficiently large $i$ and $\sum_i M^n_i < 2^{-n}$. Taking the sum of all $M^n$,
we get a lower semicomputable semimeasure $M$ such that $r_i/M_i\to 0$. Then $r_i/m_i\to 0$ also for the universal semimeasure $m$.\qed

\medskip
This result provides yet another proof that a sum of two non-random lower semicomputable reals is non-random (since the sum of two sequences that converge to~$0$ also converges to~$0$).

It shows also that Solovay functions exist (which is not immediately obvious from the definition). Moreover, it shows that there exist computable \emph{non-decreasing} Solovay functions: take a computable series of rational numbers with random sum and make this series non-increasing not changing the sum (by splitting too big terms into small pieces).

It also implies that slow convergence (in the Solovay sense) is not a property of a series itself, but only of its sum. It looks strange: some property of a computable series (of positive rational numbers), saying that \emph{infinitely many
terms come close to the upper bound provided by a priori probability}, depends only on the sum of this series. At first it seems that by splitting the terms into small parts we can destroy the property not changing the sum, but it is not so. In the next section we try to understand this phenomenon providing a direct proof for it (and as a byproduct we get some improvements in the result of this section).

\section{The Solovay property as a property of the sum}

First, let us note that the Solovay property is invariant under computable permutations. Indeed, computable permutation $\pi$ changes the a priori probability only by a constant factor: $m_{\pi(i)}= \Theta(m_i)$.
Then let us consider grouping. Since we want to allow infinite groups, let us consider a computable series $\sum_{i,j} a_{ij}$ of non-negative rational numbers. Then
	$$
\alpha=\sum_{i,j} a_{ij}=(a_{00}+a_{01}+\ldots)+(a_{10}+a_{11}+\ldots)+\ldots=\sum_i A_i,
	$$
where $A_i=\sum_j a_{ij}$.

We want to show that $A_i$ and $a_{ij}$ are slowly converging series (in the Solovay sense) at the same time. Note that slow convergence is permutation-invariant,  so it is well defined for two-dimensional series.

However, some clarifications and restrictions are needed. First, $\sum A_i$ is not in general a computable series, it is only a lower semicomputable one. We can extend the definition of the Solovay property to lower semicomputable series, still requiring $A_i=O(m_i)$, and asking this bound to be $O(1)$-tight infinitely often.  Second, such a general statement is not true: imagine that all non-negative terms are in the first group $A_0$ and all $A_1,A_2,\ldots$ are zeros. Then $\sum A_i$ does not have the Solovay property while $\sum a_{ij}$ could have it.

The following result is essentially in~\cite{HolzlKM2009}:

\begin{theorem}
	\label{Solovay-grouping}
Assume that each group $A_i$ contains only finitely many non-zero terms. Then the properties
$A_i/m_i\to 0$ and $a_{ij}/m_{ij}\to 0$ are equivalent.
\end{theorem}

Here $m_{ij}$ is the a priori probability of pair $\langle i,j\rangle$ (or its number in some computable numbering, this does not matter up to $O(1)$-factor). The convergence means that for every $\varepsilon>0$ the inequality $a_{ij}/m_{ij}>\varepsilon$ is true only for finitely many pairs $\langle i,j\rangle$.

\smallskip
\textbf{Proof}. Let us recall first that $m_i = \sum_j m_{ij}$ up to a $O(1)$-factor. (Indeed, the sum in the right hand side is lower semicomputable, so it is $O(m_i)$ due to the maximality. On the other hand, already the first term $m_{i0}$ is $\Omega(m_i)$.) So if $a_{ij}/m_{ij}$ tends to zero, the ratio $A_i/\sum_j m_{ij}$ does the same (only finitely many pairs have $a_{ij}>\varepsilon m_{ij}$ and they appear only in finitely many groups).

It remains to show that $A_i/m_i\to 0$ implies $a_{ij}/m_{ij}\to 0$. Here we need to use that only finitely many terms in each group are non-zero. For this it is enough to construct some lower semicomputable $\tilde{m}_{ij}$ such that $a_{ij}/\tilde{m}_{ij}\to 0$, somehow using the fact that $A_i/m_i\to 0$. The natural idea would be to split $m_i$ between $\tilde{m}_{ij}$ in the same proportion as $A_i$ is split between $a_{ij}$. However, for this we need to know how many terms among $a_{i0}, a_{i1},\ldots$ are non-zero, and in general this is a non-computable information. (For the special case of finite grouping this argument would indeed work.)

So we go in the other direction. For some constant $c$ we may let $\tilde{m}_{ij}$ to be $c a_{ij}$ while this does not violate the property $\sum_j \tilde{m}_{ij} \le m_i$. (When $m_i$ increases, we increase $\tilde{m}_{ij}$ when possible.) If indeed $A_i/m_i\to 0$, for every constant $c$ we have $cA_i \le m_i$ for all sufficiently large $i$, so $a_{ij}/\tilde{m}_{ij}\le 1/c$ for all sufficiently large $i$ (and only finitely many pairs $\langle i,j\rangle$ violate this requirement, because each $A_i$ has only finitely many non-zero terms). So we are close to our goal ($a_{ij}/\tilde{m}_{ij}\to 0$): it remains to perform this construction for all $c=2^{2n}$ and combine the resulting $\tilde{m}$'s with coefficients $2^{-n}$.\qed

As a corollary of Theorem~\ref{Solovay-grouping} we see (in an alternative way) that the Solovay property depends only on the sum of the series. Indeed, if $\sum_i a_i = \sum_j b_j$, these two series could be obtained by a different grouping of terms in some third series $\sum_k c_k$. To construct $c_k$, we draw intervals of lengths $a_1,a_2,\ldots$ starting from zero point, as well as the intervals of lengths $b_1,b_2\ldots$; combined endpoints split the line into intervals of lengths $c_1,c_2,\ldots$ (as shown):

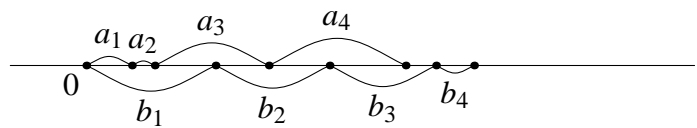
\begin{center}
       \includegraphics[scale=1.0]{lowerrandom-3.mps}
\end{center}

In this way we get not only the alternative invariance proof, but also can strengthen Theorem~\ref{Solovay-property-and-randomness}. It dealt with computable series of rational numbers. Now we still consider series of rational numbers but the summands are presented as lower semicomputable numbers and each has only finitely many different approximations. (So $r_i=\lim_n r(i,n)$ where $r$ is a computable function of $i$ and $n$ with rational values which is non-decreasing as a function of $n$ and for every $i$ there are only finitely many different values $r(i,n)$.)

Now the result of~\cite{HolzlKM2009} follows easily:

\begin{theorem}
 	\label{Solovay-property-and-randomness-strong}
Let $\alpha=\sum_i r_i$ be a converging semicomputable series of rational numbers in the sense explained above. The number $\alpha$ is random if and only if this series converges slowly in the Solovay sense \textup(i.e., $r_i/m_i$ does not converge to~$0$\textup).
\end{theorem}

\textbf{Proof}. Indeed, each $r_i$ is a sum of a computable series of non-negative rational numbers with only finitely many non-zero terms. So we can split $\sum r_i$ into a double series not changing the sum (evidently) and the Solovay property (due to Theorem~\ref{Solovay-grouping}).\qed

\medskip
In particular, we get the following corollary: \emph{an upper semicomputable function $n\mapsto f(n)$ with integer values is an upper bound for $\KP(n)$ if and only if $\sum_n 2^{-f(n)}$ is finite\textup; this
bound is tight infinitely often if and only if this sum is random}.

\medskip

Now we can show an alternative proof that all complete reals have the Solovay property.  First we observe that the Solovay property is upward closed with respect to Solovay reducibility. Indeed, if $\sum a_i$ and $\sum b_i$ are computable series of non-negative rational numbers, and $a_i$ converges slowly, then $\sum (a_i+b_i)$ converges slowly, too (its terms are bigger). So it remains to prove directly that at least one slowly converging series (or, in other terms, computable Solovay function) exists. To construct it, we watch how the values of a priori probability increase (it is convenient again to consider a priori probability of pairs):
     $$
\begin{array}{ccccc}
m_{00} & m_{01} & m_{02} & m_{03} & \ldots\\
m_{10} & m_{11} & m_{12} & m_{13} & \ldots\\
m_{20} & m_{21} & m_{22} & m_{23} & \ldots\\
\ldots & \ldots & \ldots & \ldots & \ldots
\end{array}
    $$
and fill a similar table with rational numbers $a_{ij}$ in such a way that $a_{ij}/m_{ij}\not\to0$. How do we fill this table? For each row we compute the sum of current values $m_{i,*}$; if it crosses one of the thresholds $1/2, 1/4, 1/8\ldots$, we put the crossed threshold value into the $a$-table (filling it with zeros from left to right while waiting for the next threshold crossed). In this way we guarantee that $a_{ij}$ is a computable function of $i$ and $j$; the sum of $a$-values is at most twice bigger than the sum of $m$-values; finally, in every row there exists at least one $a$-value that is at least half of the corresponding $m$-value. Logarithms of
$a$-values form a Solovay function (and $a_{ij}$ itself form a slowly convergent series).

Note that this construction does not give a \emph{nondecreasing} Solovay function directly (it seems that we still need to use the arguments from the preceding section).

\section{Busy beavers and convergence regulators}

We had several definitions that formalize the intuitive idea of a ``slowly converging series''.  However, the following one (probably the most straightforward) was not considered yet. If $a_n\to \alpha$, for every $\varepsilon>0$ there exists some $N$ such that $|\alpha-a_n|<\varepsilon$ for all $n>N$. The minimal $N$ with this property (considered as a function of $\varepsilon$, denoted by $\varepsilon\mapsto N(\varepsilon)$) is called \emph{modulus of convergence}. A sequence (or a series) should be considered ``slowly converging'' if this function grows fast. Indeed, slow convergence (defined as the Solovay property) could be equivalently characterized in these terms (see Theorem~\ref{two-slow-convergence} below).

First we define a prefix-free version of busy beaver function:

\begin{definition}
Let $m$ be a natural number. Define $BP(m)$ as the minimal value of $N$ such that $\KP(n)>m$ for all $n>N$.
\end{definition}

In other terms, $BP(m)$ is the maximal number $n$ whose prefix complexity $\KP(n)$ does not exceed $m$. Let us recall a well-known natural interpretation of $BP(m)$ in terms of ``busy beavers'':

\begin{theorem}
	Fix an optimal prefix-free universal machine $M$. Let $T(m)$ be the maximal time needed for termination of \textup(terminating\textup) programs of length at most $m$. Then
	$$
BP(m-c) \le T(m) \le BP(m+c)
	$$
for some $c$ and all $m$.
\end{theorem}

\textbf{Proof}. First we prove that for all $t>T(m)$ the compexity of $t$ is at least $m-O(1)$, thus showing that $T(m)\ge BP(m-c)$. Indeed, let $\KP(t)=m-d$. Appending the shortest program for $t$ to the prefix-free description of $d$, we get a prefix free description of the pair $\langle t,m\rangle$. Indeed, we can reconstruct $t$ and $m-d$ from the shortest program of $t$ (the second is its length) and then add $d$ and get $m$. 
Then, knowing $t$ and $m$, we run $t$ steps of all programs of length at most $m$, and then choose the first string that is not among their outputs. This string has by construction prefix complexity greater than $m$, and it is (prefix-freely) described by $m-d+O(\log d)$ bits, so $d=O(1)$.

On the other hand, $T(m)$ can be (prefix-freely) described by most long-playing program of size at most $m$ (program determines its execution time), so $\KP(T(m))\le m+O(1)$ and therefore $T(m)\le BP(m+O(1))$.\qed

\medskip
Now we can prove the equivalence of two notions of ``slow convergence'':

\begin{theorem}\label{two-slow-convergence}
The computable series of non-negative rational numbers $\sum r_i$ has the Solovay property $(\Leftrightarrow \text{has a random sum})$ if and only its modulus of convergence satisfies the inequality $N(2^{-m}) > BP(m-c)$ for some $c$ and for all $m$.
\end{theorem}

\textbf{Proof}. Let $\alpha=\sum r_i = \lim a_i$, where $a_i=r_0+\ldots+r_{i-1}$.
Assume that $\alpha$ is random. We have to show that $|\alpha-a_i|<2^{-m}$ implies $\KP(i)>m-O(1)$; this shows that $N(2^{-m})\ge BP(m-O(1))$. Since $\KP(i)=\KP(a_i)+O(1)$, it is enough to show that every rational $2^{-m}$-approximation to $\alpha$ has complexity at least $m-O(1)$. This is a bit stronger condition than the condition $\KP(\alpha_0\ldots \alpha_{m-1})\ge m-O(1)$ (used in prefix complexity version of Schnorr--Levin theorem) since now we consider \textsl{all} approximations, not only the prefix of the binary expansion. However, it can be proven in a similar way.

Let $c$ be some integer.  Consider an effectively open set $U_c$ constructed as follows. For every rational $r$ we consider the neighborhood around $r$ of radius $2^{-\KP(r)-c}$; the set $U_c$ is the union of these neighborhoods.  (Since $\KP(r)$ is upper semicomputable, it is indeed an effectively open set.) The total length of all intervals is $2\cdot 2^{-c} \sum_r 2^{-\KP(r)} \le 2^{-(c-1)}$. Therefore, $U_c$ form a Martin-L\"of test, and random $\alpha$ does not belong to $U_c$ for some $c$. This means that complexity of $2^{-m}$-ap\-pro\-xi\-ma\-t\-ions of $\alpha$ is at least $m-O(1)$.

In the other direction we can use Schnorr--Levin theorem without any changes: if $N(2^{-m})\ge BP(m-c)$, then $\KP(i)\ge m-O(1)$ for every $i$ such that $a_i$ is a $2^{-m}$-app\-roxi\-mation to $\alpha$. Therefore, the $m$-bit prefix of $\alpha$ has complexity at least $m-O(1)$, since knowing this prefix we can effectively find an $a_i$ that exceeds it (and the corresponding~$i$).\qed

\medskip
\textbf{Question}. Note that this theorem shows equivalence between two formalizations of an intuitive idea of ``slowly converging series'' (or three, if we consider the Solovay reducibility as a way to compare the rate of convergence). However, the proof goes through Martin-L\"of randomness of the sum (where the series itself disappears). Can we have a more direct proof? Can we connect the Solovay reducibility (not only completeness) to the properties of the modulus of convergence?

\medskip

Reformulating the definition of $BP(m)$ in terms of a priori probability, we say that $BP(m)$ is the minimal $N$ such that all $n>N$ have a priori probability less than $2^{-m}$. However, in terms of a priori probability the other definition looks more natural: let $BP'(m)$ be the minimal $N$ such that the total a priori probability of all $n>N$ is less than $2^{-m}$. Generally speaking, $BP'(m)$ can be greater that $BP(m)$, but it turns out that it still can be used to characterize randomness in the same way:

\begin{theorem}
	\label{competing-beavers}
Let $a_i$ be a computable increasing sequence of rational numbers that converges to a random number $\alpha$. Then $N(2^{-m})\ge BP'(m-c)$.
\end{theorem}

\textbf{Proof}. Since all $i>N(2^{-m})$ have the same a priori probability as the corresponding $a_i$ (up to some $O(1)$-factor), it is enough to show that for every $m$ the sum of a priori probabilities of all rational numbers in the $2^{-m}$-neighborhood of a random~$\alpha$ is $O(2^{-m})$ (recall that for all $i>N(2^{-m})$ the corresponding $a_i$ belong to this neighborhood).

As usual, we go in the other direction and cover all ``bad'' $\alpha$ that do not have this property by a set of small measure.  Not having this property means that for every $c$ there exists $m$ such that the sum of a priori probabilities of rational numbers in the $2^{-m}$-neighborhood of $\alpha$ exceeds $c2^{-m}$. For a given $c$, we consider all intervals with rational endpoints that have the following property: \emph{the sum of a priori probabilities of all rational numbers in this interval is more than $c/2$ times bigger than the interval's length}. Every bad $\alpha$ is covered by an interval with this property (the endpoints of the interval $(\alpha-2^{-m},\alpha+2^{-m})$ can be changed slightly to make them rational), and the set of intervals having this property is enumerable. It is enough to show that the union of all such intervals has measure $O(1/c)$, in fact, at most $4/c$.

It is also enough to consider a finite union of intervals with this property. Moreover, we may assume that this union does not contain redundant intervals (that can be deleted without changing the union). Let us order all the intervals according to their left endpoints:
	$$
(l_0,r_0), (l_1,r_1),(l_2,r_2),\ldots
	$$
where $l_0\le l_1\le l_2\le\ldots$\;
It is easy to see that right endpoints go in the same order (otherwise one of the intervals would be redundant). So $r_0\le r_1\le r_2\le\ldots$ Now note that $r_i\le l_{i+2}$, otherwise the interval $(l_{i+1},r_{i+1})$ would be redundant. Therefore, intervals with even numbers $(l_0,r_0), (l_2,r_2), (l_4,r_4)\ldots$ are disjoint, and for each of them the length is $c/2$ times less than the sum of a priori probabilities of rational numbers inside it. Therefore, the total length of these intervals does not exceed $2/c$, since the sum of all priori probabilities is at most~$1$. The same is true for intervals with odd numbers, so in total we get the bound~$4/c$.\qed

\medskip
\textbf{Question}: We see that both $BP$ and $BP'$ can be used to characterize randomness, but how much could $BP$ and $BP'$ differ in general?

\bigskip

\textbf{Acknowledgments.} The authors are grateful to the organizers of the conference for the opportunity to present this work, and to L.~Staiger and R.~H\"olzl for comments.

\bibliographystyle{plain}	
\bibliography{lowerrandom}

\end{document}